\documentclass[10pt,letterpaper]{article}

\makeatletter

\usepackage{latexsym}

\makeatother

\usepackage{blindtext,titlefoot}
\usepackage{authblk}
\usepackage{marvosym}
\usepackage[utf8]{inputenc}
\usepackage[english]{babel}
\usepackage{amsmath}
\usepackage{amsfonts}
\usepackage{amssymb}
\usepackage{makeidx}
\usepackage{graphicx}
\usepackage{lmodern}
\usepackage{kpfonts}
\usepackage{dsfont}
\usepackage{cancel}
\usepackage{amsthm} 
\usepackage[left=2cm,right=2cm,top=2cm,bottom=2cm]{geometry}
\usepackage{subcaption}
\usepackage{multirow}
\usepackage{float}
\usepackage{url}
\usepackage{breakurl}
\usepackage[breaklinks]{hyperref}
\usepackage{multicol}
\usepackage{mathtools, cuted}
\usepackage{lipsum}
\usepackage{booktabs}
\usepackage{comment}
\usepackage{cite}

\usepackage{fancyhdr}
\fancypagestyle{firstpage}{\fancyhf{}
\setcounter{page}{1}
\rfoot{\thepage}
}

\providecommand{\nset}[1]{
\mathbb{#1}
}
\providecommand{\set}[1]{
\left\{#1\right\}
}

\providecommand{\ifr}[5]{
{}^{#1}_{#2}{#3}_{#4}^{#5}
}
\providecommand{\gam}[1]{
\Gamma\left(#1 \right)
}
\providecommand{\re}[1]{
Re\left(#1 \right)
}
\providecommand{\norm}[1]{
\left\lVert #1 \right\rVert
}
\providecommand{\abs}[1]{
\left\lvert #1 \right\rvert
}
\providecommand{\ds}[1]{
\displaystyle #1
}
\providecommand{\der}[3]{
\dfrac{#1^{#3} }{ #1 #2^{#3}}
}

\providecommand{\re}[1]{
\hbox{Re}\left(#1 \right)
}
\providecommand{\im}[1]{
\hbox{Im}\left(#1 \right)
}
\providecommand{\rnd}[2]{
\hbox{Rnd}_#2\left(#1\right)
}

\newtheorem{theorem}{ Theorem}[section]
\newtheorem{definition}[theorem]{Definition}
\newtheorem{proposition}[theorem]{Proposition}
\newtheorem{corollary}[theorem]{Corollary}
\newtheorem{example}[theorem]{Example}

\usepackage{enumitem} 
\setlist[itemize]{noitemsep} 

\usepackage{titlesec} 
\titleformat{\section}[block]{\large\bfseries\scshape\centering}{\thesection.}{1em}{} 
\titleformat{\subsection}[block]{\large\bfseries\scshape\centering}{\thesubsection.}{1em}{}
\titleformat{\subsubsection}[block]{\large\bfseries\scshape\centering}{\thesubsubsection.}{1em}{} 

\title{\huge\bfseries Reduction of a nonlinear system and its numerical solution using a fractional iterative method}

\author[,a]{A. Torres-Hernandez  \footnote{Email address: anthony.torres@ciencias.unam.mx; ORCID: 0000-0001-6496-9505}}
\affil[a]{\normalsize Department of Physics, Faculty of Science - UNAM, Mexico}

\author[,b]{F. Brambila-Paz \footnote{Email address: fernandobrambila@gmail.com; ORCID: 0000-0001-7896-6460}}
\affil[b]{\normalsize Department of Mathematics, Faculty of Science - UNAM, Mexico}

\affil[c]{\normalsize Faculty of Engineering, Universidad Panamericana - Aguascalientes, Mexico}

\author[,c,d]{P. M. Rodrigo \footnote{Email address: prodrigo@up.edu.mx; ORCID: 0000-0003-0100-6124}}
\affil[d]{\normalsize Centre for Advanced Studies in Energy and Environment (CEAEMA), University of Jaén, Spain.}

\author[,c]{E. De-la-Vega \footnote{Email address: evega@up.edu.mx}}

\date{}

\begin{document}

\maketitle

\thispagestyle{firstpage}


\begin{abstract}

A nonlinear algebraic equation system of $5$ variables is numerically solved, which allows modeling the behavior of the temperatures and the efficiencies of a hybrid solar receiver, which in simple terms is the combination of a photovoltaic system with a thermoelectric system. In addition, a way to reduce the previous system to a nonlinear system of only $ 2 $ variables is presented. Naturally, reducing algebraic equation systems of dimension $ N $ to systems of smaller dimensions has the main advantage of reducing the number of variables involved in a problem, but the analytical expressions of the systems become more complicated. However, to minimize this disadvantage, an iterative method that does not explicitly depend on the analytical complexity of the system to be solved is used.  A fractional iterative method, valid for one and several variables, that uses the properties of fractional calculus, in particular the fact that the fractional derivatives of constants are not always zero, to find solutions of nonlinear systems is presented.

\end{abstract}

\section{Introduction}

A classic problem in mathematics, which is of common interest in physics and engineering, is finding the set of zeros of a function $f:\Omega \subset \nset{R}^n \to \nset{R}^n$, that is,

\begin{eqnarray}\label{eq:1-001}
\set{\xi \in \Omega \ : \ \norm{f(\xi)}=0},
\end{eqnarray}

where $\norm{ \ \cdot \ }: \nset{R}^n \to \nset{R}$ denotes any vector norm. Although finding the zeros of a function may seem like a simple problem, in general it involves solving an \textbf{algebraic equation system} as follows

\begin{eqnarray}\label{eq:1-002}
\left\{
\begin{array}{c}
\left[f\right]_1(x)=0\\
\left[f\right]_2(x)=0\\
\vdots \\
\left[f\right]_n(x)=0
\end{array}\right.,
\end{eqnarray}

where $[f]_k: \nset{R}^n \to \nset{R}$ denotes the $k$-th component of the function $f$.  Let $\set{\hat{e}_i}_{i=1}^n$ be the canonical basis of $\nset{R}^n$, if the nature of function $f$ allows it, it is possible to find a function $g_1:\nset{R}^n \to \nset{R}^n$ that allows rewriting the system of equations \eqref{eq:1-002} as follows

\begin{eqnarray}\label{eq:1-003}
\left\{
\begin{array}{c}
\left[x\right]_1- [g_1]_1\left( x-\hat{e}_1[x]_1 \right)=0\\
\left[x\right]_2- [g_1]_2\left( x-\hat{e}_2[x]_2 \right)=0\\
\vdots \\
\left[x\right]_n- [g_1]_n\left( x-\hat{e}_n[x]_n \right)=0
\end{array}\right.,
\end{eqnarray}

it should be noted that the system of equations \eqref{eq:1-003} may represent a \textbf{linear system} or a \textbf{nonlinear system}. Without loss of generality, taking the variable $[x]_1$ and replacing it in the rest of equations, the previous system may be rewritten as

\begin{eqnarray}\label{eq:1-0031}
\left\{
\begin{array}{c}
\left[x\right]_2- [g_1]_2\left( x-\hat{e}_2[x]_2-\hat{e}_1[x]_1+\hat{e}_1 [g_1]_1\left( x-\hat{e}_1[x]_1 \right) \right)=0\\
\left[x\right]_3- [g_1]_3\left( x-\hat{e}_3[x]_3-\hat{e}_1[x]_1+\hat{e}_1 [g_1]_1\left( x-\hat{e}_1[x]_1 \right) \right)=0\\
\vdots \\
\left[x\right]_n- [g_1]_n\left( x-\hat{e}_n[x]_n-\hat{e}_1[x]_1+\hat{e}_1 [g_1]_1\left( x-\hat{e}_1[x]_1 \right) \right)=0
\end{array}\right.,
\end{eqnarray}

if the nature of the function $g_1$ allows it, it is possible to find a function $g_2:\nset{R}^n \to \nset{R}^n$, with $[g_2]_1 \equiv 0$, that allows rewriting the system of equations \eqref{eq:1-0031} as follows

\begin{eqnarray}\label{eq:1-0032}
\left\{
\begin{array}{c}
\left[x\right]_2- [g_2]_2\left( x-\hat{e}_2[x]_2-\hat{e}_1[x]_1 \right)=0\\
\left[x\right]_3- [g_2]_3\left( x-\hat{e}_3[x]_3-\hat{e}_1[x]_1 \right)=0\\
\vdots \\
\left[x\right]_n- [g_2]_n\left( x-\hat{e}_n[x]_n -\hat{e}_1[x]_1\right)=0
\end{array}\right.,
\end{eqnarray}

which corresponds to an equivalent system of the system of equations \eqref{eq:1-003}. Assuming that the process to build the systems \eqref{eq:1-0031} and \eqref{eq:1-0032} may be repeated successively up to the variable $[x]_{k-2}$, it is possible to obtain the following system of equations

\begin{eqnarray}\label{eq:1-0033}
\left\{
\begin{array}{c}
\left[x\right]_{k-1}- [g_{k-1}]_{k-1}\left( x-\hat{e}_{k-1}[x]_{k-1}- \sum_{r=1}^{k-2}\hat{e}_{r}[x]_{r} \right)=0\\
\left[x\right]_k- [g_{k-1}]_k\left( x-\hat{e}_k[x]_k - \sum_{r=1}^{k-2}\hat{e}_{r}[x]_{r}\right)=0\\
\vdots \\
\left[x\right]_n- [g_{k-1}]_n\left( x-\hat{e}_n[x]_n- \sum_{r=1}^{k-2}\hat{e}_{r}[x]_{r} \right)=0
\end{array}\right.,
\end{eqnarray}

if by taking the variable $[x]_{k-1}$ and replacing it in the rest of equations, it is not possible to obtain an equivalent system of the system of equations \eqref{eq:1-003}, it is possible to find a  function $h:\nset{R}^n \to \nset{R}^n$, with $[h]_r \equiv 0 \ \forall r<k $, that allows rewriting the resulting system as follows

\begin{eqnarray}\label{eq:1-004}
\left\{
\begin{array}{c}
\left[x\right]_{k}- [h]_{k}\left( x- \sum_{r=1}^{k-1}\hat{e}_{r}[x]_{r} \right)=0\\
\left[x\right]_{k+1}- [h]_{k+1}\left( x- \sum_{r=1}^{k-1}\hat{e}_{r}[x]_{r} \right)=0\\
\vdots \\
\left[x\right]_n- [h]_n\left( x- \sum_{r=1}^{k-1}\hat{e}_{r}[x]_{r} \right)=0
\end{array}\right.,
\end{eqnarray}

the system of equations \eqref{eq:1-004} represents a \textbf{transcendental system}, that is, there are no algebraic operations that allow rewriting \eqref{eq:1-004} in an equivalent system of the system of equations \eqref{eq:1-003}.  It should be noted that the previous system has the advantage of needing fewer variables, although its analytical expression becomes more complicated. However, it is possible to extract the same information as that contained in the system \eqref{eq:1-002}. In general, it is necessary to use numerical methods of the iterative type to approximate to the solution of the system \eqref{eq:1-004}.

It is necessary to mention that the iterative methods have an intrinsic problem, since if a system has $N$ solutions it is necessary to invest time in finding $N$ initial conditions, but this problem is partially solved by combining iterative methods with fractional calculus, whose result is known as \textbf{fractional iterative methods}, because these new methods have the ability to find $N$ solutions of a system using a single initial condition. In this document, a fractional iterative method that does not explicitly depend on the analytical complexity or the fractional partial derivatives of the function for which zeros are searched is presented, then it is an ideal iterative method for solving nonlinear systems in several variables.

\section{Fixed Point Method}

Let  $\Phi:\nset{R}^n \to \nset{R}^n$ be a function. It is possible to build a sequence $\set{x_i}_{i=0}^\infty$  by defining the following iterative method

\begin{eqnarray}\label{eq:2-001}
x_{i+1}:=\Phi(x_i),
\end{eqnarray}

if is true that $x_i\to \xi\in \nset{R}^n$ and if the function $\Phi$ is continuous around $\xi$, we obtain that

\begin{eqnarray}\label{eq:2-002}
\xi=\lim_{i\to \infty}x_{i+1}=\lim_{i\to \infty}\Phi(x_i)=\Phi\left(\lim_{i\to \infty}x_i \right)=\Phi(\xi),
\end{eqnarray}

the above result is the reason by which the method \eqref{eq:2-001} is known as the \textbf{fixed point method}. Moreover, the function $\Phi$ is called an \textbf{iteration function}. To understand the nature of the convergence of the iteration function $\Phi$, the following definition is necessary \cite{plato2003concise}

\begin{definition}
Let $\Phi:\nset{R}^n \to \nset{R}^n$  be an iteration function. The method \eqref{eq:2-001} for determining $\xi\in \nset{R}^n$ is called \textbf{(locally) convergent}, if there exists $\delta>0$ such that for every initial value

\begin{eqnarray*}
x_0\in B(\xi;\delta):=\set{y\in \nset{R}^n \ : \  \norm{y-\xi}<\delta},
\end{eqnarray*}

it holds that

\begin{eqnarray}\label{eq:2-003}
\lim_{i \to \infty}\norm{x_i-\xi}\to 0 & \Rightarrow & \lim_{i\to \infty}x_i=\xi.
\end{eqnarray}

\end{definition}

If we have a function $f:\Omega \subset \nset{R}^n \to \nset{R}^n$ for which we want to determine the set \eqref{eq:1-001}, in general it is possible to write an iteration function 
$\Phi$ as follows \cite{burden2002analisis}

\begin{eqnarray*}
\Phi(x)=x-A(x)f(x),
\end{eqnarray*}

with $A(x)$ a matrix given as follows

\begin{eqnarray*}
A(x)=:\left([A]_{jk}(x) \right)=\begin{pmatrix}
[A]_{11}(x)&[A]_{12}(x)& \cdots &[A]_{1n}(x)\\
[A]_{21}(x)&[A]_{22}(x)&\cdots& [A]_{2n}(x)\\
\vdots & \vdots & \ddots & \vdots \\
[A]_{n1}(x)&[A]_{n2}(x)&\cdots &[A]_{nn}(x)
\end{pmatrix},
\end{eqnarray*}

where $[A]_{jk}(x):\nset{R}^n \to \nset{R}$. It is necessary to mention that the matrix  $ A (x) $ is determined according to the order of convergence desired \cite{torreshern2020}.

\subsection{Order of Convergence}

Consider the following definition \cite{plato2003concise}

\begin{definition}
Let $ \Phi: \Omega \subset \nset{R}^ n \to \nset{R}^ n $ be an iteration function with a fixed point $ \xi \in \Omega $. Then the method \eqref{eq:2-001} is called  \textbf{(locally) convergent of (at least) order $ \boldsymbol{p} $} ($ p \geq 1 $), if there are exists $ \delta> 0 $  and $ C $, a non-negative constant  with $ C <1 $ if $ p = 1 $, such that for any initial value $ x_0 \in B (\xi; \delta) $ it holds that

\begin{eqnarray}\label{eq:c2.08}
\norm{x_{k+1}-\xi}\leq C \norm{x_k-\xi}^p, & k=0,1,2,\cdots,
\end{eqnarray}

where $ C $ is called convergence factor.

\end{definition}

The order of convergence is usually related to the speed at which the sequence generated by \eqref{eq:2-001} converges. For the particular case $ p = 1 $ it is said that the method \eqref{eq:2-001} has an \textbf{order of convergence (at least) linear}, and for the case $p=2$  it is said that the method \eqref{eq:2-001} has an \textbf{order of convergence (at least) quadratic}. The following theorem, allows characterizing the order of convergence of an iteration function $ \Phi $ with its derivatives \cite{plato2003concise,stoer2013}. Before continuing,  we need to consider the following multi-index notation. Let $\nset{N}_0$ be the set $\nset{N}\cup\set{0}$, if $\gamma \in \nset{N}_0^n$ then

\begin{eqnarray}
\left\{
\begin{array}{l}
 \gamma!:= \ds\prod_{k=1}^n [\gamma]_k ! \vspace{0.1cm}\\
 \abs{\gamma}:= \ds \sum_{k=1}^n [\gamma]_k\vspace{0.1cm}\\
 x^\gamma:= \ds \prod_{k=1}^n [x]_k^{[\gamma]_k}\vspace{0.1cm}\\
\der{\partial}{x}{\gamma}:= \dfrac{\partial^{\abs{\gamma}}}{\partial [x]_1^{[\gamma]_1}\partial [x]_2^{[\gamma]_2}\cdots \partial [x]_n^{[\gamma]_n} }
\end{array}\right. .
\end{eqnarray}

\begin{theorem}\label{teo:c2.01}
Let $ \Phi: \Omega \subset \nset{R}^n \to \nset {R}^n $ be an iteration function with a fixed point $ \xi \in \Omega $. Assuming that $\Phi $ is $ p$-times differentiable in $ \xi $ for some $ p \in \nset{N} $, and moreover

\begin{eqnarray}\label{eq:c2.09}
\left\{
\begin{array}{cc}
\ds  \dfrac{\partial^\gamma [\Phi]_k(\xi) }{ \partial x^\gamma}=0, \ \forall k\geq 1 \mbox{ and } \forall  \abs{\gamma}<p, & \mbox{if }p\geq 2 \vspace{0.1cm}\\
\ds \norm{\Phi^{(1)}(\xi)}<1, & \mbox{if }p=1
\end{array}\right.,
\end{eqnarray}

where $\Phi^{(1)}$ denotes the \textbf{Jacobian matrix} of the function $\Phi$, then $ \Phi $ is (locally) convergent of (at least) order $ p $.

\begin{proof}

Let $\Phi:\nset{R}^n \to \nset{R}^n$ be an iteration function, and let $\set{\hat{e}_k}_{k=1}^n$ be the canonical basis of $\nset{R}^n$. Considering the following index notation (Einstein notation)

\begin{eqnarray*}
\Phi(x)=\sum_{k=1}^n [\Phi]_k(x)\hat{e}_k: = [\Phi]_k(x)\hat{e}_k=\hat{e}_k[\Phi]_k(x),
\end{eqnarray*}

and using the Taylor series expansion of a vector-valued function in multi-index notation, we obtain two cases:

\begin{itemize}

\item[i)]  Case $p\geq 2:$

\begin{align*}
\Phi(x_i)
=& \ds  \Phi(\xi)+  \sum_{\abs{\gamma} =1}^p \dfrac{1}{\gamma !}\hat{e}_k\dfrac{\partial^\gamma [\Phi]_k(\xi) }{ \partial x^\gamma}   (x_i-\xi)^\gamma   + \hat{e}_k[o]_k\left(\max_{\abs{\gamma}=p}  \set{ (x_i-\xi)^\gamma }\right) \\
=& \ds \Phi(\xi)+ \sum_{m =1}^p   \left( \sum_{\abs{\gamma} =m}\dfrac{1}{\gamma !} \hat{e}_k \dfrac{\partial^\gamma [\Phi]_k(\xi) }{ \partial x^\gamma}   (x_i-\xi)^\gamma \right)   + \hat{e}_k[o]_k\left(\max_{\abs{\gamma}=p}  \set{ (x_i-\xi)^\gamma }\right),
\end{align*}

then

\begin{align*}
\norm{\Phi(x_i)-\Phi(\xi)}&\leq  \ds \sum_{m =1}^p   \left( \sum_{\abs{\gamma} =m}\dfrac{1}{\gamma !}\norm{ \hat{e}_k \dfrac{\partial^\gamma [\Phi]_k(\xi) }{ \partial x^\gamma}   (x_i-\xi)^\gamma }\right)  +  \norm{\hat{e}_k  [o]_k\left(\max_{\abs{\gamma}=p}  \set{ (x_i-\xi)^\gamma }\right)} \\
&\leq   \ds \sum_{m =1}^p   \left( \sum_{\abs{\gamma} =m}\dfrac{1}{\gamma !}\norm{ \dfrac{\partial^\gamma [\Phi]_k(\xi) }{ \partial x^\gamma}  \hat{e}_k   }\right)\norm{x_i-\xi}^m+  o\left( \norm{x_i-\xi}^p \right), 
\end{align*}

assuming that $\xi$ is a fixed point of $\Phi$ and that $\dfrac{\partial^\gamma [\Phi]_k(\xi) }{ \partial x^\gamma}=0 \ \forall k\geq 1$ and $  \forall \abs{\gamma}<p$ is fulfilled, the previous expression implies that

\begin{eqnarray*}
\dfrac{\norm{\Phi(x_i)-\Phi(\xi)}}{\norm{x_i-\xi}^p}=\dfrac{\norm{x_{i+1}-\xi}}{\norm{x_i-\xi}^p}\leq\sum_{\abs{\gamma} =p}\dfrac{1}{\gamma !}\norm{ \dfrac{\partial^\gamma [\Phi]_k(\xi) }{ \partial x^\gamma}  \hat{e}_k   } +\dfrac{o\left(\norm{x_i-\xi}^p \right)}{\norm{x_i-\xi}^p},
\end{eqnarray*}

therefore

\begin{eqnarray*}
\lim_{i\to \infty} \dfrac{\norm{x_{i+1}-\xi}}{\norm{x_i-\xi}^p}\leq  \sum_{\abs{\gamma} =p}\dfrac{1}{\gamma !}\norm{ \dfrac{\partial^\gamma [\Phi]_k(\xi) }{ \partial x^\gamma}  \hat{e}_k   },
\end{eqnarray*}

as a consequence, if the sequence $\set{x_i}_{i=0}^\infty$ generated by \eqref{eq:2-001} converges to $\xi$, there exists a value $k>0$ such that

\begin{eqnarray*}
\norm{x_{i+1}-\xi}\leq \left( \sum_{\abs{\gamma} =p}\dfrac{1}{\gamma !}\norm{ \dfrac{\partial^\gamma [\Phi]_k(\xi) }{ \partial x^\gamma}  \hat{e}_k   }\right)\norm{x_i-\xi}^p, & \forall i\geq k,
\end{eqnarray*}

then $ \Phi $ is (locally) convergent of (at least) order $ p $.

\item[ii)] Case $p=1:$

\begin{align*}
\Phi(x_i)
=& \ds  \Phi(\xi)+   \sum_{\abs{\gamma} =1} \dfrac{1}{\gamma !}\hat{e}_k \dfrac{\partial^\gamma [\Phi]_k(\xi) }{ \partial x^\gamma}   (x_i-\xi)^\gamma  +  \hat{e}_k  [o]_k\left(\max_{\abs{\gamma}=1}  \set{ (x_i-\xi)^\gamma } \right) \\
=& \ds  \Phi(\xi)+ \Phi^{(1)}(\xi) (x_i-\xi) +  \hat{e}_k[o]_k\left(\max_{\abs{\gamma}=1}  \set{ (x_i-\xi)^\gamma }\right)  ,
\end{align*}

then

\begin{align*}
\norm{\Phi(x_i)-\Phi(\xi)}\leq  \norm{\Phi^{(1)}(\xi)} \norm{x_i-\xi }+   o\left( \norm{x_i-\xi} \right) ,
\end{align*}

assuming that $\xi$ is a fixed point of $\Phi$ and that $\norm{\Phi^{(1)}(\xi)}<1$ is fulfilled, the previous expression implies that

\begin{eqnarray*}
\dfrac{\norm{\Phi(x_i)-\Phi(\xi)}}{\norm{x_i-\xi}}=\dfrac{\norm{x_{i+1}-\xi}}{\norm{x_i-\xi}}\leq\norm{\Phi^{(1)}(\xi)} +\dfrac{o\left(\norm{x_i-\xi} \right)}{\norm{x_i-\xi}},
\end{eqnarray*}

therefore

\begin{eqnarray*}
\lim_{i\to \infty} \dfrac{\norm{x_{i+1}-\xi}}{\norm{x_i-\xi}}\leq\norm{\Phi^{(1)}(\xi)},
\end{eqnarray*}

as a consequence, if the sequence $\set{x_i}_{i=0}^\infty$ generated by \eqref{eq:2-001} converges to $\xi$, there exists a value $k>0$ such that

\begin{eqnarray*}
\norm{x_{i+1}-\xi}\leq \norm{\Phi^{(1)}(\xi)}\norm{x_i-\xi}, & \forall i\geq k,
\end{eqnarray*}

then $ \Phi $ is (locally) convergent of order (at least) linear.

\end{itemize}

\end{proof}

\end{theorem}

The following proposition follows from the previous theorem

\begin{proposition}
Let $f:\Omega \subset \nset{R}^n\to \nset{R}^n$ be a function with a value $\xi\in \Omega$ such that $\norm{f(\xi)}=0$, and let $\Phi:\nset{R}^n \to \nset{R}^n$ be an iteration function as follows

\begin{eqnarray}\label{eq:2-004}
\Phi(x)=x-A(x)f(x),
\end{eqnarray}

with $A(x)$ a matrix. If the following condition is fulfilled

\begin{eqnarray}\label{eq:c2.17}
\lim_{x\to \xi}A(x)=\left(f^{(1)}(\xi) \right)^{-1},
\end{eqnarray}

then $\Phi$ satisfies a necessary (but not sufficient) condition to be (locally) convergent of order (at least) quadratic in $B(\xi;\delta)$.

\begin{proof}
From the \textbf{Theorem \ref{teo:c2.01}} we have that an iteration function has an order of convergence (at least) quadratic if it fulfills the following condition

\begin{eqnarray*}
\lim_{x\to \xi}\dfrac{\partial [\Phi]_k(x)}{\partial [x]_j}=0, & \forall j,k\leq n,
\end{eqnarray*}

which may be written equivalently as follows

\begin{eqnarray}\label{eq:c2.161}
 \lim_{x\to \xi}\norm{\Phi^{(1)}(x)}=0.
\end{eqnarray}

Then, we can assume that we have a function $f(x):\Omega \subset \nset{R}^n \to \nset{R}^n$  with a zero $\xi \in \Omega$, such that all of its first partial derivatives are defined in $ \xi $, and taking the iteration function $ \Phi $ given by \eqref{eq:2-004}, the $ k $-th component of the iteration function may be written as

\begin{eqnarray*}
[\Phi]_k(x)=[x]_k-\sum_{j=1}^n[A]_{kj}(x)[f]_j(x),
\end{eqnarray*}

then

\begin{eqnarray*}
\partial_l [\Phi]_k(x)=\delta_{lk}-\sum_{j=1}^n \big{(} [A]_{kj}(x)\partial_l [f]_j(x)+\left(\partial_l [A]_{kj}(x) \right)[f]_j(x) \big{)},
\end{eqnarray*}

where $ \delta_{lk} $ is the Kronecker delta, which is defined as

\begin{eqnarray*}
\delta_{lk}=\left\{
\begin{array}{cc}
1,& \mbox{si }l=k\\
0,& \mbox{si }l\neq k
\end{array}\right..
\end{eqnarray*}

Assuming that \eqref{eq:c2.161} is fulfilled

\begin{eqnarray*}
\partial_l [\Phi]_k(\xi)=\delta_{lk}-\sum_{j=1}^n [A]_{kj}(\xi)\partial_l [f]_j(\xi)=0 & \Rightarrow & \sum_{j=1}^n [A]_{kj}(\xi)\partial_l [f]_j(\xi)=\delta_{lk}, \ \forall l,k\leq n,
\end{eqnarray*}

the previous expression may be written in matrix form as

\begin{eqnarray*}
A(\xi) f^{(1)}(\xi)=I_n & \Rightarrow & A(\xi)= \left(f^{(1)}(\xi)\right)^{-1},
\end{eqnarray*}

where  $ I_n $ denotes the identity matrix of $ n \times n $. Then any matrix $ A (x) $ that fulfills the following condition

\begin{eqnarray*}
\lim_{x\to \xi}A(x)= \left(f^{(1)}(\xi)\right)^{-1}
\end{eqnarray*}

guarantees that exists $ \delta> 0 $ such that iteration function $ \Phi $ given by \eqref{eq:2-004} satisfies a necessary (but not sufficient) condition to be (locally) convergent of  order (at least) quadratic in $ B (\xi; \delta) $.
\end{proof}

\end{proposition}

Finally, the following corollary follows from the previous proposition

\begin{corollary}\label{cor:2-001}
Let $\Phi:\nset{R}^n \to \nset{R}^n$ be an iteration function. If $\Phi$ defines a sequence $\set{x_i}_{i=0}^\infty$ such that $x_i\to \xi$, and if the following condition is fulfilled

\begin{eqnarray}\label{eq:c2.16}
\lim_{x\to \xi}\norm{\Phi^{(1)}(x)}\neq 0,
\end{eqnarray}

then $\Phi$ has an order of convergence (at least) linear  in $B(\xi;\delta)$.
\end{corollary}

\section{Riemann-Liouville Fractional Derivative}

One of the key pieces in the study of fractional calculus is the iterated integral, which is defined as follows \cite{hilfer00}

\begin{definition}
Let $ L_{loc} ^ 1 (a, b) $ be the space of locally integrable functions in the interval $ (a, b) $. If $ f $ is a function such that $ f \in L_ {loc} ^ 1 (a, \infty) $, then the $n$-th iterated integral of the function $ f $ is given by 

\begin{eqnarray}\label{eq:c1.16}
\begin{array}{c}
\ds \ifr{}{a}{I}{x}{n} f(x)=\ifr{}{a}{I}{x}{}\left(\ifr{}{a}{I}{x}{n-1} f(x)  \right)=\frac{1}{(n-1)!}\int_a^x(x-t)^{n-1}f(t)dt,
\end{array}
\end{eqnarray}

where

\begin{eqnarray*}
\ifr{}{a}{I}{x}{} f(x):=\int_a^x f(t)dt.
\end{eqnarray*}

\end{definition}

Considerate that $ (n-1)! = \gam{n} $
, a generalization of \eqref{eq:c1.16} may be obtained for an arbitrary order $ \alpha> 0 $

\begin{eqnarray}\label{eq:c1.17}
\ifr{}{a}{I}{x}{\alpha} f(x)=\dfrac{1}{\gam{\alpha}}\int_a^x(x-t)^{\alpha-1}f(t)dt,
\end{eqnarray}

the equation \eqref{eq:c1.17}  correspond to the definition of \textbf{Riemann-Liouville (right) fractional integral}. Fractional integrals satisfy the  \textbf{semigroup property}, which is given in the following proposition \cite{hilfer00}

\begin{proposition}
Let $ f $ be a function. If $ f \in L_{loc} ^ 1 (a, \infty) $, then the fractional integrals of $ f $ satisfy that

\begin{eqnarray}\label{eq:c1.19}
\ifr{}{a}{I}{x}{\alpha} \ifr{}{a}{I}{x}{\beta}f(x) = \ifr{}{a}{I}{x}{\alpha + \beta}f(x),& \alpha,\beta>0.
\end{eqnarray}

\end{proposition}

From the previous proposition, and considering that the operator $ d / dx $  is the inverse operator to the left of the operator $ \ifr {}{a}{I}{x}{} $, any integral $ \alpha$-th of a function $ f \in L_{loc} ^ 1 (a, \infty) $ may be written as

\begin{eqnarray}\label{eq:c1.20}
\ifr{}{a}{I}{x}{\alpha}f(x)=\dfrac{d^n}{dx^n}\ifr{}{a}{I}{x}{n}\left( \ifr{}{a}{I}{x}{\alpha}f(x) \right)=\dfrac{d^n}{dx^n}\left( \ifr{}{a}{I}{x}{n+\alpha}f(x)\right).
\end{eqnarray}

With the previous results, we can build the operator  \textbf{Riemann-Liouville fractional derivative} as follows \cite{hilfer00,kilbas2006theory}

\begin{eqnarray}\label{eq:c1.23}
\normalsize
\begin{array}{c}
\ifr{}{a}{D}{x}{\alpha}f(x) := \left\{
\begin{array}{cc}
\ds \ifr{}{a}{I}{x}{-\alpha}f(x), &\mbox{if }\alpha<0 \vspace{0.1cm}\\  
\ds \dfrac{d^n}{dx^n}\left( \ifr{}{a}{I}{x}{n-\alpha}f(x)\right), & \mbox{if }\alpha\geq 0
\end{array}
\right.
\end{array}, 
\end{eqnarray}

where  $ n = \lfloor \alpha \rfloor + 1 $. Considering $a=0$, then applying the  operator \eqref{eq:c1.23} to the  function $ x^{\mu} $, with  $\alpha \in \nset{R}\setminus\nset {Z} $ and $\mu>-1$, we obtain the following result

\begin{eqnarray}\label{eq:c1.13}
\ifr{}{0}{D}{x}{\alpha}x^\mu = 
 \dfrac{\gam{\mu+1}}{\gam{\mu-\alpha+1}}x^{\mu-\alpha}.
\end{eqnarray}

\section{Fractional Pseudo-Newton Method}

Let $f:\Omega \subset \nset{R}^n\to \nset{R}^n$ be a function. We can consider the problem of finding a value $\xi\in \Omega$ such that $\norm{f(\xi)}$=0. A first approximation to value $\xi$ is by a linear approximation of the function $f$ in a valor $x_i\in \Omega$ with $\norm{x_i- \xi} < \epsilon$, that is,

\begin{eqnarray}\label{eq:c2.36}
f(x)\approx f(x_i)+f^{(1)}(x_i)(x-x_i),
\end{eqnarray}

then considering that $ \xi $ is a zero of $ f $, from the previous expression we obtain that

\begin{eqnarray*}
0\approx f(x_i)+f^{(1)}(x_i)(\xi-x_i) & \Rightarrow &  \xi \approx x_i- \left(f^{(1)}(x_i) \right)^{-1} f(x_i),
\end{eqnarray*}

consequently, may be generated a sequence $ \set{x_i}_{i = 0} ^ \infty $ that approximates the value $ \xi $  using the iterative method

\begin{eqnarray*}
x_{i+1}:=\Phi(x_i)=x_i- \left(f^{(1)}(x_i) \right)^{-1} f(x_i), & i=0,1,2,\cdots,
\end{eqnarray*}

which corresponds to well-known Newton's method. However, the equation \eqref{eq:c2.36} is not the only way to generate a linear approximation to the  function $ f $ in the point $ x_i $, another alternative is to use the next approximation

\begin{eqnarray}\label{eq:c2.37}
f(x)\approx f(x_i)+mI_n(x-x_i),
\end{eqnarray}

where $ I_n $ corresponds to the identity matrix of $n\times n$ and $ m $ is any constant value of a slope, that allows the approximation  \eqref{eq:c2.37} to the  function $ f $ to be valid. The previous equation allows to obtain the following iterative method

\begin{eqnarray}\label{eq:c2.39}
x_{i+1}:=\Phi(x_i)= x_i- \left( m^{-1}I_n \right) f(x_i), & i=0,1,2\cdots,
\end{eqnarray}

which corresponds to a particular case of the \textbf{parallel chord method} \cite{ortega1970iterative}. It is necessary to mention that for some definitions of fractional derivative, it is fulfilled that the derivative of the order $ \alpha $ of a constant is different from zero, that is,

\begin{eqnarray}\label{eq:c2.30}
\partial_k^\alpha c :=\der{\partial}{[x]_k}{\alpha}c \neq 0 , & c=constant,
\end{eqnarray}

where $ \partial_k ^ \alpha $ denotes any fractional derivative applied only in the component $ k $, that does not cancel the constants and that fulfills the following continuity relation with respect to the order $ \alpha $ of the derivative

\begin{eqnarray}\label{eq:c2.301}
\lim_{\alpha \to 1}\partial_k^\alpha c=\partial_kc.
\end{eqnarray}

Considering a function $\Phi:(\nset{R}\setminus \nset{Z})\times \nset{C}^n \to \nset{C}^n$. Then, using as a basis the idea of the method \eqref{eq:c2.39}, and considering any fractional derivative that fulfills the conditions \eqref{eq:c2.30} and \eqref{eq:c2.301}, we can define the \textbf{fractional pseudo-Newton method} as follows

\begin{eqnarray}\label{eq:c2.401}
x_{i+1}:=\Phi(\alpha, x_i)= x_i- P_{\epsilon,\beta}(x_i) f(x_i), & i=0,1,2\cdots,
\end{eqnarray}

with $\alpha\in\nset{R}\setminus\nset{Z}$, in particular $\alpha\in[-2,2]\setminus\nset{Z}$ \cite{torreshern2020}, where $ P_{\epsilon, \beta} (x_i) $ is a matrix evaluated in the value $ x_i $, which is given by the following expression

\begin{eqnarray}\label{eq:c2.402}
P_{\epsilon,\beta}(x_i):=\left([P_{\epsilon,\beta}]_{jk}(x_i)\right)=\left( \partial_k^{\beta(\alpha,[x_i]_k)}\delta_{jk}+ \epsilon\delta_{jk}  \right)_{x_i},
\end{eqnarray}

where

\begin{eqnarray}
\partial_k^{\beta(\alpha,[x_i]_k)}\delta_{jk}:= \der{\partial}{[x]_k}{\beta(\alpha,[x_i]_k)}\delta_{jk}, & 1\leq j,k\leq n,
\end{eqnarray}

with $ \delta_{jk} $ the Kronecker delta, $ \epsilon $ a positive constant $ \ll 1 $, and $ \beta (\alpha, [x_i]_k) $ a function defined as follows

\begin{eqnarray}\label{eq:c2.34}
\beta(\alpha,[x_i]_k):=\left\{
\begin{array}{cc}
\alpha, &\mbox{if \hspace{0.1cm} }  |[x_i]_k|\neq 0 \vspace{0.1cm}\\
1,& \mbox{if \hspace{0.1cm}  }  |[x_i]_k|=0
\end{array}\right. .
\end{eqnarray}

Due to the part of the integral operator that fractional derivatives usually have, we consider in the matrix \eqref{eq:c2.402} that each fractional derivative is obtained for a real variable $[x]_k$, and if the result allows it, this variable is subsequently substituted by a complex variable $[x_i]_k$. It should be mentioned that the value $ \alpha = 1 $ in \eqref{eq:c2.34}, is taken to avoid the discontinuity that is generated when using the fractional derivative of constants in the value $ x = 0 $. Moreover, since in the previous method  $\norm{\Phi^{(1)}(\alpha,\xi)}\neq 0$ if $\norm{f(\xi)}=0$, for the \textbf{Corollary \ref{cor:2-001}}, any sequence $ \set{x_i} _ {i = 0} ^ \infty $ generated by the iterative method \eqref {eq:c2.401} has an order of convergence (at least) linear.

To finish this section, it is necessary to mention that although the interest in fractional calculus has mainly focused on the study and development of techniques to solve differential equation systems of order non-integer \cite{hilfer00,kilbas2006theory, martinez2017applications1,martinez2017applications2,torreshern2019proposal}. Over the years, iterative methods have also been developed that use the properties of fractional derivatives to solve algebraic equation systems \cite{gao2009local,torres2017fractional,brambila2018fractional,akgul2019fractional,torreshern2020,torres2020fractional,gdawiec2020newton,torres2020approximation}. 
These methods may be called \textbf{fractional iterative methods}, which under certain conditions, may accelerate their speed of convergence with the implementation of the Aitken's method     \cite{stoer2013,brambila2018fractional}. 

It should be noted that depending on the definition of fractional derivative used, fractional iterative methods have the particularity that they may be used of local form \cite{gao2009local} or of global form \cite{torreshern2020}. These methods also have the peculiarity of being able to find complex roots of polynomials using real initial conditions \cite{torres2017fractional}.  Some differences between Newton's method and two fractional iterative methods are listed in the Table \ref{tab:07}

\begin{table}[!ht]
\centering
\footnotesize
\begin{tabular}{c|c|c|c}
\toprule
&Classical Newton&Fractional Newton& Fractional Pseudo-Newton \\ \midrule
\begin{tabular}{c}
Can it find complex zeros\\
 of a polynomial using \\
 real initial conditions?
\end{tabular}& No & Yes & Yes \\ \midrule
\begin{tabular}{c}
Can it find multiple zeros \\
of a function using a \\
single initial condition?
\end{tabular}& No & Yes & Yes \\ \midrule
\begin{tabular}{c}
Can it be used if the function \\ 
is not differentiable?
\end{tabular}& No&Yes&Yes \\ \midrule
\begin{tabular}{c}
For a space of dimension $ N $ \\
are needed
\end{tabular}&
\begin{tabular}{c}
$ N \times N $ classic \\
 partial derivatives
\end{tabular}& 
\begin{tabular}{c}
$ N \times N $ fractional \\
 partial derivatives
\end{tabular}
 & \begin{tabular}{c}
$ N$ fractional \\
 partial derivatives
\end{tabular} \\ \midrule
\begin{tabular}{c}
Is it recommended for solving systems\\
 where the (fractional) partial derivatives  \\
 are analytically difficult to obtain?
\end{tabular}& No & No & Yes\\ \bottomrule
\end{tabular}
\caption{Some differences between the classical Newton's method and two fractional iterative methods.}\label{tab:07}
\end{table}

\subsection{Some Examples}

Instructions for implementing the method \eqref{eq:c2.401} along with information to provide values $\alpha \in [-2,2]\setminus \nset{Z}$ are found in the reference \cite{torreshern2020}. For rounding reasons, for the examples the following function is defined

\begin{eqnarray}\label{eq:4-001}
\rnd{[x]_k}{m}:=\left\{
\begin{array}{cc}
\re{[x]_k},& \mbox{ if \hspace{0.1cm}} \abs{\im{[x]_k}}\leq 10^{-m}\vspace{0.1cm}\\
\left[x\right]_k,& \mbox{ if \hspace{0.1cm}} \abs{\im{[x]_k}}> 10^{-m}\vspace{0.1cm}\\
\end{array}\right..
\end{eqnarray}

Combining the function \eqref{eq:4-001} with the method \eqref{eq:c2.401}, the following iterative method is defined

\begin{eqnarray}\label{eq:c2.40}
x_{i+1}:=\rnd{\Phi(\alpha, x_i)}{5}, & i=0,1,2\cdots.
\end{eqnarray}

\begin{example}

Let $\set{f_k}_{k=0}^\infty$ be a sequence of functions, with

\begin{eqnarray*}
f_k(x)=\dfrac{\pi}{2}-\sum_{m=0}^k \dfrac{(-1)^m x^{2 m+1}}{(2 m+1) \Gamma(2 m+2)} & \underset{k \to \infty}{\longrightarrow} & \int_x^\infty \dfrac{\sin(t)}{t}dt.
\end{eqnarray*}

Then considering the value $k=50$, the initial condition $x_0=1.85$ is chosen to use the iterative method  given by \eqref{eq:c2.40} along with fractional derivative given by \eqref{eq:c1.13}. Consequently, we obtain the results of the Table \ref{tab:01}

\begin{table}[!ht]
\centering
\footnotesize
$
\begin{array}{c|ccccc}
\toprule
&\alpha& x_n&\norm{x_n - x_{n-1} }_2  &\norm{f_{50}\left(x_n \right)}_2& n \\ 
\midrule
1	&	-0.83718	&	23.60399266        	&	4.10000e-7	&	9.80551e-9	&	30	\\
 2	&	 -0.81526	&	  29.87824476        	&	7.30000e-7	&	8.19133e-7	&	 273	\\
3 	&	 -0.71339	&	  17.33566366        	&	4.10000e-7	&	2.48591e-8	&	  34	\\
 4	&	 -0.71324	&	  11.08303768        	&	3.70000e-7	&	2.67499e-8	&	  29	\\
 5	&	 -0.71174	&	   4.89383571        	&	5.10000e-7	&	4.87621e-8	&	  24	\\
 6	&	  0.36251	&	 -12.29964074 - 4.38965942i	&	4.41814e-7	&	9.42697e-7	&	  42	\\
 7	&	  0.36333	&	 -31.27978791 - 5.29112884i	&	4.31045e-7	&	5.52647e-7	&	 391	\\
 8	&	  0.36684	&	 -24.97153098 - 5.07020771i	&	3.56931e-7	&	7.71425e-7	&	  42	\\
 9	&	  0.38451	&	 -24.97153097 + 5.07020788i	&	2.05913e-7	&	5.36982e-7	&	  44	\\
 10 	&	  0.38646	&	 -18.65002028 + 4.78651268i	&	2.30000e-7	&	5.54454e-7	&	  42	\\
11	&	  0.44711	&	  -5.86005858 - 3.72373544i	&	4.72017e-7	&	9.90371e-7	&	  41	\\
12	&	  0.55885	&	 -31.27978639 + 5.2911368i 	&	5.50000e-7	&	3.69789e-7	&	 183	\\
13	&	  1.41172	&	   1.92644561  	&	1.10000e-7	&	9.97696e-7	&	 196	\\
\bottomrule
\end{array}
$
\caption{Results obtained using the iterative method \eqref{eq:c2.40} with $\epsilon=e-3$.}\label{tab:01}
\end{table}

\end{example}

\begin{example}

Let $f$ be a function, with

\begin{eqnarray*}
f(x)=
\begin{pmatrix}
\dfrac{1}{2} [x]_1\big{(} \sin\left( [x]_1 [x]_2 \big{)} -1\right)-\dfrac{1}{4\pi}[x]_2 \vspace{0.1cm}\\
\left( 1- \dfrac{1}{4\pi} \right)\left(e^{2[x]_1}-e\right) +e\left(\dfrac{1}{\pi}[x]_2-2[x]_1\right)
\end{pmatrix}.
\end{eqnarray*}

Then the initial condition $x_0=(0.86,0.86)^T$ is chosen to use the iterative method  given by \eqref{eq:c2.40} along with fractional derivative given by \eqref{eq:c1.13}. Consequently, we obtain the results of the Table \ref{tab:02}

\begin{table}[!ht]
\centering
\footnotesize
$
\begin{array}{c|cccccc}
\toprule
&\alpha& [x_n]_1& [x_n]_2 &\norm{x_n - x_{n-1} }_2  &\norm{f\left(x_n \right)}_2& n \\ 
\midrule
1 	&	 0.7283 	&	 -0.13780202 - 0.87180273i	&	   2.16460988 - 4.68221226i	&	9.11043e-8	&	8.81449e-7	&	 100	\\
 2	&	 0.72889	&	 -0.15442216        	&	   1.14021866        	&	6.22977e-7	&	8.30511e-7	&	  60	\\
 3	&	 0.78188	&	 -0.20477864 - 1.30850366i	&	   2.21623485 - 7.86783099i	&	5.56776e-8	&	9.92736e-7	&	 246	\\
 4	&	 0.86097	&	  1.14584377 + 0.68994256i	&	   8.09450017 - 5.99607116i	&	2.64575e-8	&	9.42041e-7	&	 249	\\
 5	&	 1.11159	&	  1.70987637        	&	 -18.87534307        	&	1.41421e-8	&	9.92487e-7	&	 447	\\
 6	&	 1.14766	&	  1.48216448        	&	  -8.41311536        	&	1.41421e-8	&	8.86632e-7	&	 233	\\
 7 	&	 1.17262	&	 -1.36674692 + 0.07786741i	&	     -5.76423 + 0.47853094i	&	2.00000e-8	&	9.92337e-7	&	 394	\\
8	&	 1.18538	&	 -1.36674698 - 0.07786726i	&	  -5.76422966 - 0.4785315i 	&	2.23607e-8	&	9.88600e-7	&	 387	\\
9	&	 1.19954	&	  1.57643706        	&	   -12.098725        	&	1.41421e-8	&	7.09538e-7	&	 386	\\
10	&	 1.20058	&	  1.64946521        	&	 -15.55495398        	&	1.41421e-8	&	9.10544e-7	&	 465	\\
11	&	 1.2852 	&	 -0.76073057 + 0.14192444i	&	  -2.11123992 + 0.82667655i	&	1.02470e-7	&	8.39720e-7	&	  97	\\
12	&	 1.29642	&	  1.34362303    	&	  -4.29400761      	&	7.61577e-8	&	4.60872e-7	&	  92	\\
\bottomrule
\end{array}
$
\caption{Results obtained using the iterative method \eqref{eq:c2.40} with $\epsilon=e-3$.}\label{tab:02}
\end{table}

\end{example}

\begin{example}

Let $f$ be a function, with

\begin{eqnarray*}
f(x)=
\begin{pmatrix}
-3.6[x]_2\left(\cos\left( [x]_2^2 \right)+[x]_1^3[x]_3 \right) - 3.6[x]_3 + 10.8 \vspace{0.1cm}\\
-1.6[x]_1\left([x]_1 + [x]_2^3[x]_3\right) - 1.6\sinh \left([x]_3\right) + 6.4 \vspace{0.1cm}\\
-4.6[x]_2\left( [x]_1[x]_3^3 + 1\right) - 4.6\cosh \left([x]_1 \right) + 27.6
\end{pmatrix}.
\end{eqnarray*}

Then the initial condition $x_0=(0.95,0.95,0.95)^T$ is chosen to use the iterative method  given by \eqref{eq:c2.40} along with fractional derivative given by \eqref{eq:c1.13}. Consequently, we obtain the results of the Table \ref{tab:03}

\begin{table}[!ht]
\centering
\footnotesize
$
\begin{array}{c|cccccccc}
\toprule
&\alpha& [x_n]_1& [x_n]_2& [x_n]_3 &\norm{x_n - x_{n-1} }_2  &\norm{f\left(x_n \right)}_2& n \\ 
\midrule
1	&	0.96828	&	-0.28991424 + 1.38566039i	&	0.4041105 - 1.39254282i	&	-0.62409681 + 1.25568859i	&	2.44949e-8	&	8.34134e-7	&	78	\\
2 	&	 0.9698 	&	  0.62792492 - 1.29495978i	&	  0.57678001 - 1.30742987i	&	  0.48322895 - 1.29731024i	&	2.64575e-8	&	9.12958e-7	&	  64	\\
 3	&	 0.96985	&	 -0.28991423 - 1.38566037i	&	  0.40411047 + 1.39254282i	&	 -0.62409683 - 1.25568861i	&	2.00000e-8	&	9.32021e-7	&	  83	\\
 4	&	 0.97106	&	 -0.58277447 + 0.49660576i	&	 -0.49995795 + 1.39319334i	&	  0.09221108 - 1.69837571i	&	3.00000e-8	&	6.23838e-7	&	  62	\\
 5	&	 0.97192	&	  0.62792485 + 1.29495977i	&	  0.57677997 + 1.3074299i 	&	  0.48322899 + 1.29731024i	&	3.31662e-8	&	9.61719e-7	&	  57	\\
 6	&	 0.97823	&	 -0.12415396 + 0.98083552i	&	 -0.51004547 - 1.39105393i	&	 -0.57743861 - 1.50487453i	&	2.44949e-8	&	9.15509e-7	&	 270	\\
 7	&	 0.97858	&	 -0.12415386 - 0.98083557i	&	 -0.51004543 + 1.3910539i 	&	 -0.57743856 + 1.50487453i	&	1.41421e-8	&	9.86401e-7	&	 278	\\
 8	&	 1.03775	&	  1.30219735        	&	 -1.31677799        	&	  -1.4605226        	&	2.82843e-8	&	8.76439e-7	&	  98	\\
 9	&	 1.04019	&	 -1.43433659        	&	  1.27415875        	&	 -1.11130559        	&	4.89898e-8	&	8.39041e-7	&	  54	\\
 10 	&	 1.0421 	&	 -1.16248344 + 0.0469604i 	&	 -0.62570099 - 0.42962177i	&	  1.74938849 - 0.27012065i	&	3.31662e-8	&	8.92465e-7	&	 279	\\
11	&	 1.96396	&	  0.53848559 - 0.36927367i	&	  0.64776248 + 0.48376485i	&	  2.00930932 - 0.07078346i	&	1.73205e-8	&	9.76607e-7	&	 216	\\
12	&	 1.96537	&	   0.5384856 + 0.3692736i 	&	  0.64776247 - 0.48376478i	&	  2.00930935 + 0.07078343i	&	1.41421e-8	&	9.37835e-7	&	 208	\\
\bottomrule
\end{array}
$
\caption{Results obtained using the iterative method \eqref{eq:c2.40} with $\epsilon=e-3$.}\label{tab:03}
\end{table}

\end{example}

\section{Equations of a Hybrid Solar Receiver}

Considering the notation

\begin{eqnarray*}
(T_{cell},T_{hot},T_{cold},\eta_{cell},\eta_{TEG})^T:=\left([x]_1,[x]_2,[x]_3,[x]_4,[x]_5 \right)^T,
\end{eqnarray*}

it is possible to define the following system of equations that corresponds to the combination of a solar photovoltaic system with a thermoelectric generator system \cite{bjork2015performance,bjork2018maximum}, which is named as a \textbf{hybrid solar receiver}

\begin{eqnarray}\label{eq:004}
\left\{
\begin{array}{l}
\left[x\right]_1=[x]_2+a_1\cdot a_2\left( 1-[x]_4 \right)\\
\left[x\right]_2=[x]_3+a_1\cdot a_3 \left( 1-[x]_4\right)\left(1-[x]_5\right)\\
\left[x\right]_3=a_4+a_1\cdot a_5 \left( 1-[x]_4\right)\left(1-[x]_5\right)\\
\left[x\right]_4 =a_6[x]_1+a_7\\
\left[x\right]_5=(a_8-1)\left(1-\dfrac{[x]_3+a_9}{[x]_2+a_9} \right)\left(a_8+ \dfrac{[x]_3+a_9}{[x]_2+a_9}\right)^{-1}
\end{array}\right.,
\end{eqnarray}

whose deduction and some details about the difficulty in finding its solution may be found in the reference \cite{rodrigo2019performance}. The $ a_i $'s in the previous system are constants defined by the following expressions

\begin{eqnarray*}
\left\{
\begin{array}{l}
\begin{array}{l}
 a_2=r_{cell}+r_{sol}+A_{cell}\left(\dfrac{r_{cop}+r_{cer}}{A_{TEG}}+\dfrac{r_{intercon}}{0.5\cdot \sqrt{f^*\cdot A_{TEG}}\left(b\cdot \sqrt{f^*}+\sqrt{A_{TEG}} \right) } \right)\\
a_5=A_{cell}\left( \dfrac{r_{intercon}}{0.5\cdot \sqrt{f^*\cdot A_{TEG}}\left(b\cdot \sqrt{f^*}+\sqrt{A_{TEG}} \right) }+\dfrac{r_{cer}}{A_{TEG}}+R_{heat\_ exch} \right)
\end{array}\\
\begin{array}{lll}
a_1=\eta_{opt}\cdot C_g \cdot DNI,&
a_3=\dfrac{A_{cell}\cdot l}{f^*\cdot A_{TEG}\cdot k_{TEG}},&a_4=T_{air} \\ \\
a_6=-\eta_{cell,ref}\cdot \gamma_{cell},&a_7=\eta_{cell,ref}\left(1+25 \cdot \gamma_{cell} \right),&
a_8=\sqrt{1+ZT}\\ \\
a_9=273.15&
\end{array}
\end{array}\right.,
\end{eqnarray*}

with the following particular values \cite{rodrigo2019performance}

\begin{eqnarray*}
\left\{
\begin{array}{lll}
\eta_{opt}=0.85, & r_{intercon}=2.331 e-7 ,& T_{air}=20 \\
        C_g=800, &     A_{cell}=9e-6, &     R_{heat\_exch}=0.5 \\
DNI=900 ,&     A_{TEG}=5.04e-5, &     \eta_{cell,ref}=0.43 \\
r_{cell}=3e-6, &     f^*=0.7, &     \gamma_{cell}=4.6e-4 \\
r_{sol}=1.603e-6, & b=5e-4, &     ZT=1 \\
r_{cop}=7.5e-7, &     l=5e-4 ,& r_{cer}=8e-6\\
k_{TEG}=1.5 
\end{array}
\right..
\end{eqnarray*}

Using the system of equations \eqref{eq:004}, it is possible to define a  function $f_1:\Omega \subset \nset{R}^5\to \nset{R}^5$, that is,

\begin{eqnarray}\label{eq:005}
f_1(x):=\begin{pmatrix}
\left[x\right]_1-[x]_2-a_1\cdot a_2\left( 1-[x]_4 \right)\\
\left[x\right]_2-[x]_3-a_1\cdot a_3 \left( 1-[x]_4\right)\left(1-[x]_5\right)\\
\left[x\right]_3-a_4-a_1\cdot a_5 \left( 1-[x]_4\right)\left(1-[x]_5\right)\\
\left[x\right]_4 -a_6[x]_1-a_7\\
\left[x\right]_5-(a_8-1)\left(1-\dfrac{[x]_3+a_9}{[x]_2+a_9} \right)\left(a_8+ \dfrac{[x]_3+a_9}{[x]_2+a_9}\right)^{-1}
\end{pmatrix}.
\end{eqnarray}

Then the initial condition $x_0=(53.67,51.82,21.54,0.43,0.01)^T$ is chosen to use the iterative method  given by \eqref{eq:c2.40} along with fractional derivative given by \eqref{eq:c1.13}. Consequently, we obtain the results of the Table \ref{tab:04}

\begin{table}[!ht]
\centering
\footnotesize
$
\begin{array}{c|ccccccccc}
\toprule
&\alpha&[x_n]_1&[x_n]_2&[x_n]_3&[x_n]_4&[x_n]_5&||x_n-x_{n-1}||_2&||f_1(x_n)||_2&n \\ \midrule
1&1.02632&53.76229916&51.55509481&22.07807195&0.42431082&0.01618411 	&9.35752e-6	&9.99890e-3	&4112
 \\ \bottomrule
\end{array}
$
\caption{Results obtained using the iterative method \eqref{eq:c2.40} with $\epsilon=e-4$.}\label{tab:04}
\end{table}

\subsection{Reducing the Number of Equations of a Hybrid Solar Receiver}

It should be noted that system \eqref{eq:004} corresponds to a particular case of the system of equations \eqref{eq:1-003}, as a consequence through consecutive substitutions of certain variables $[x]_k$'s may be transformed into a particular case of the system of equations \eqref{eq:1-004}. In particular, if we do consecutive substitutions of variables $[x]_1, \ [x]_4, \ [x]_5$ and some algebraic simplifications, we may obtain the following transcendental system

\begin{eqnarray}\label{eq:006}
\left\{
\begin{array}{l}
\left[x\right]_2=[x]_3-a_1\cdot a_3 \dfrac{\left( a_6 [x]_2 + a_7 - 1 \right) \left( a_8 \left([x]_3+ a_9 \right) + \left([x]_2+ a_9\right)  \right)  }{(1+a_1 a_2 a_6 ) \left( a_8 \left( [x]_2 + a_9  \right) + \left([x]_3+ a_9\right) \right) } \vspace{0.1cm} \\
\left[x \right]_3=a_4-a_1\cdot a_5 \dfrac{\left( a_6 [x]_2 + a_7 - 1 \right) \left( a_8 \left([x]_3+ a_9 \right) + \left([x]_2+ a_9\right)  \right)  }{(1+a_1 a_2 a_6 ) \left( a_8 \left( [x]_2 + a_9  \right) + \left([x]_3+ a_9\right) \right) } 
\end{array}\right.,
\end{eqnarray}

whose solution allows to know the values of the variables $[x]_1,[x]_4$ and $[x]_5$ through the following equations

\begin{eqnarray}\label{eq:0061}
\left\{
\begin{array}{l}
\left[x\right]_1=\dfrac{[x]_2 - a_1 a_2 (a_7 - 1)}{1+a_1 a_2 a_6} \vspace{0.1cm}\\
\left[ x \right]_4=\dfrac{a_6 \left( a_1 a_2 + [x]_2\right) + a_7}{1+ a_1 a_2 a_6 } \vspace{0.1cm}\\
\left[ x \right]_5=\dfrac{(a_8 - 1)\left( [x]_2 - [x]_3 \right)}{a_8 \left([x]_2+a_9\right) + \left( [x]_3+a_9 \right)}
\end{array}\right..
\end{eqnarray}

Using the system of equations \eqref{eq:006}, it is possible to define a  function $f_2:\Omega \subset \nset{R}^2\to \nset{R}^2$, that is,

\begin{eqnarray}\label{eq:007}
f_2(x):=\begin{pmatrix}
\left[x\right]_2-[x]_3+a_1\cdot a_3 \dfrac{\left( a_6 [x]_2 + a_7 - 1 \right) \left( a_8 \left([x]_3+ a_9 \right) + \left([x]_2+ a_9\right)  \right)  }{(1+a_1 a_2 a_6 ) \left( a_8 \left( [x]_2 + a_9  \right) + \left([x]_3+ a_9\right) \right) } \vspace{0.1cm} \\
\left[x \right]_3-a_4+a_1\cdot a_5 \dfrac{\left( a_6 [x]_2 + a_7 - 1 \right) \left( a_8 \left([x]_3+ a_9 \right) + \left([x]_2+ a_9\right)  \right)  }{(1+a_1 a_2 a_6 ) \left( a_8 \left( [x]_2 + a_9  \right) + \left([x]_3+ a_9\right) \right) } 
\end{pmatrix}.
\end{eqnarray}

Then the initial condition $x_0=(53,19)^T$ is chosen to use the iterative method  given by \eqref{eq:c2.40} along with fractional derivative given by \eqref{eq:c1.13}. Consequently, we obtain the results of the Table \ref{tab:05}

\begin{table}[!ht]
\centering
\footnotesize
$
\begin{array}{c|cccccc}
\toprule
&\alpha&[x_n]_2&[x_n]_3&||x_n-x_{n-1}||_2&||f_2(x_n)||_2&n \\ \midrule
1&1.17778&51.55653453 &22.0782978 &3.72194e-5&	9.97906e-3&	1420
 \\ \bottomrule
\end{array}
$
\caption{Results obtained using the iterative method \eqref{eq:c2.40} with $\epsilon=e-4$.}\label{tab:05}
\end{table}

Finally, using \eqref{eq:0061} and the results of the Table \ref{tab:05}, we obtain the following values

\begin{eqnarray}\label{eq:7-001}
\left\{
\begin{array}{l}
\left[x_n\right]_1=53.76173931 \vspace{0.1cm}\\
\left[ x_n \right]_4=0.42431093 \vspace{0.1cm}\\
\left[ x_n \right]_5=0.01618472
\end{array}\right..
\end{eqnarray}

It should be mentioned that the previous results are practically the same as those presented in Table \ref{tab:04}, which shows that it is possible to extract the information contained in system \eqref{eq:004} using the system \eqref{eq:006}.

The system of equations \eqref{eq:004} generally depends on two parameters, the direct normal irradiance ($DNI$) and the ambient temperature ($T_{air}$). These parameters are measured in real-time at certain times of the day \cite{rodrigo2019performance} and it is necessary to calculate a new system solution for each new pair of parameters, that is,

\begin{eqnarray*}
(DNI,T_{air})\overset{f_1}{\longrightarrow} x_n \in \nset{R}^5.
\end{eqnarray*}

Unfortunately this system presents a certain degree of instability \cite{torres2020fractional}, as a consequence, it is difficult to find a solution through iterative methods unless a suitable initial condition is used. The latter implies that before trying to solve the system \eqref{eq:004} using some iterative method, the task of searching for a suitable initial condition in a space of $5$ dimensions must be carried out, and this task complicates obtaining in real-time the behavior of the temperatures and the efficiencies of the hybrid solar receiver. Although reducing the system dimensions does not eliminate the task of searching for an initial condition, it is less complicated to search for an initial condition in a space of 2 dimensions than in a space of 5 dimensions. Therefore, the solutions of the system \eqref{eq:006}  may be determined more quickly, and the behavior in real-time of the temperatures and the efficiencies of the hybrid solar receiver may be obtained with greater precision.

We finish this section by presenting real measurements of parameters $DNI$ and $T_{air}$, as well as the solutions they generate for the system \eqref{eq:006} using the iterative method  given by \eqref{eq:c2.40} along with fractional derivative given by \eqref{eq:c1.13}.

\begin{table}[!ht]
\centering
\footnotesize
$
\begin{array}{c|cccc|cccccc}
\toprule
&DNI&T_{air}&[x_0]_2&[x_0]_3&\alpha&[x_n]_2&[x_n]_3&||x_n-x_{n-1}||_2&||f_2(x_n)||_2&n \\ \midrule
    1     & 574.319 & 16.832 & 33    & 22    & 1.15517 & 36.97006552 & 18.16597807 & 2.14742e-5 & 9.99513e-3 & 1697 \\
    2     & 81.348 & 23.332 & 32    & 21    & 1.22759 & 26.1945392 & 23.51244959 & 3.78783e-5 & 9.98177e-3 & 1423 \\
    3     & 421.637 & 17.061 & 35    & 14    & 1.17586 & 31.86785373 & 18.03494236 & 3.41155e-5 & 9.98047e-3 & 1264 \\
    4     & 370.62 & 15.34 & 25    & 19    & 1.2069 & 28.35083227 & 16.20876627 & 5.35705e-5 & 9.96884e-3 & 997 \\
    5     & 63.796 & 19.527 & 30    & 22    & 1.23793 & 21.78372255 & 19.67591543 & 4.24110e-5 & 9.96582e-3 & 1861 \\
    6     & 173.964 & 13.955 & 23    & 12    & 1.22759 & 20.06382043 & 14.35167375 & 6.20276e-5 & 9.99870e-3 & 778 \\
    7     & 60.219 & 21.911 & 22    & 25    & 1.23793 & 24.04433152 & 22.05249253 & 3.74858e-5 & 9.97890e-3 & 1924 \\
    8     & 158.031 & 16.98 & 27    & 23    & 1.23793 & 22.54666358 & 17.34765403 & 4.07181e-5 & 9.99593e-3 & 1942 \\
    9     & 73.3474 & 31.87 & 30    & 33    & 1.21724 & 34.46203488 & 32.04766382 & 2.46233e-5 & 9.99229e-3 & 1960 \\
    10    & 114.473 & 19.798 & 20    & 15    & 1.23793 & 23.81769247 & 20.06413324 & 3.79117e-5 & 9.97919e-3 & 1927 \\
    11    & 337.977 & 14.505 & 27    & 12    & 1.24828 & 26.35550694 & 15.28940708 & 3.32622e-5 & 9.99058e-3 & 1599 \\
    12    & 290.003 & 17.678 & 25    & 20    & 1.2069 & 27.86466037 & 18.35620882 & 3.54168e-5 & 9.98776e-3 & 1178 \\
    13    & 198.9558 & 20.146 & 30    & 19    & 1.21724 & 27.15103929 & 20.60390244 & 3.60974e-5 & 9.96025e-3 & 1247 \\
    14    & 142.54 & 32.932 & 41    & 32    & 1.21724 & 37.96815767 & 33.25643301 & 2.40786e-5 & 9.96843e-3 & 1978 \\
    15    & 831.497 & 27.259 & 58    & 27    & 1.18621 & 56.50580249 & 29.1848537 & 2.60029e-5 & 9.99290e-3 & 1895 \\
    16    & 839.482 & 23.023 & 51    & 27    & 1.16552 & 52.50277786 & 24.97861392 & 2.19814e-5 & 9.98513e-3 & 1785 \\
    17    & 30.275 & 21.416 & 17    & 18    & 1.23793 & 22.47208318 & 21.48561688 & 4.08991e-5 & 9.99714e-3 & 1851 \\
    18    & 374.688 & 18.493 & 34    & 18    & 1.2069 & 31.66818022 & 19.36351722 & 2.63287e-5 & 9.98794e-3 & 1756 \\
    19    & 94.3555 & 28.373 & 35    & 28    & 1.21724 & 31.71360529 & 28.59186834 & 2.64235e-5 & 9.98198e-3 & 1978 \\
     \bottomrule
\end{array}
$
\caption{Results obtained using the iterative method \eqref{eq:c2.40} with $\epsilon=e-4$.}\label{tab:06}
\end{table}

\section{Conclusions}

The reduction of an algebraic equation system has among its advantages the fact that it is less complicated to determine initial conditions and that its solutions may be determined more quickly. However, these advantages are overshadowed by the fact that the analytical expression of the system becomes more complicated. For this reason, it is necessary to have an iterative method that is not affected by the increase in analytical complexity of the system. The fractional pseudo-Newton method does not explicitly depend on the analytical complexity or the fractional partial derivatives of the function for which zeros are searched, these characteristics make this method ideal to be implemented when reducing and solving nonlinear systems in several variables.

Partially funded by PAPPIT$\-$IT$\_$101421, UNAM.

\bibliography{Biblio}
\bibliographystyle{unsrt}

\nocite{philipps2015current,green2019solar} \nocite{perez2018efficiency} \nocite{kost2013levelized,talavera2016worldwide}
\nocite{rowe2018crc} \nocite{beeri2015hybrid,tamaki2017hybrid,kil2017highly}
\nocite{rezania2016coupled}
\nocite{munoz2015efficiencies}

\end{document}